%% file: main.tex
\documentclass{IEEEtaes}

\usepackage{color,array,amsthm}
\usepackage{graphicx}
\usepackage[T1]{fontenc}

\jvol{XX}
\jnum{XX}
\jmonth{XXXXX}
\paper{1234567}
\pubyear{2026}
\doiinfo{TAES.2022.Doi Number}

\usepackage{amsmath}
\usepackage{amssymb}
\usepackage{hyperref}

\begin{document}

\title{Time-Optimal Collision Avoidance Via a Greedy Polynomial Backward Sweep}

\author{ZENO PAVANELLO, FRANK DE VELD, and ROBERTO ARMELLIN}


\author{
	ZENO PAVANELLO\\
    Department of Aerospace Science and Technology (DAER), Politecnico di Milano, via Privata Giuseppe La Masa 34, Milano, 20156, Italy \\
    \and
    FRANK DE VELD\\
	Netherlands Aerospace Centre (NLR), Anthony Fokkerweg 2, Amsterdam, 1059 CM, The Netherlands\\
	\and 
	ROBERTO ARMELLIN\\
	Te Pu\=naha \=Atea Auckland - Space Institute (TPA-SI), The University of Auckland, 20 Symonds Street, Auckland, 1010, New Zealand     
}
\maketitle

\begin{abstract}
	Spacecraft collision avoidance for low-thrust satellites often requires determining not only how to maneuver, but also how late a maneuver can begin while still ensuring safety. This paper presents a greedy time-optimal (GTO) backward-sweep method to find the latest maneuver initiation time. The method starts from the nominal time of closest approach and iteratively propagates the maneuver backward in time, selecting at each step the thrust direction that locally minimizes the chosen danger metric. Differential algebra is used to efficiently propagate state sensitivities and update the time of closest approach online. The method is tested on a large dataset of conjunctions, using both miss distance and probability of collision as safety metrics. The approach achieves accurate results and only a small loss of optimality relative to an optimal-control benchmark, while retaining runtimes suitable for on-board implementation.
\end{abstract}
\glsresetall


\section{INTRODUCTION}
Recently, the increasing presence of space debris in Earth-centered orbits has posed significant challenges for satellite operators. As space debris orbits Earth uncontrollably, it can pose a threat of collision to nearby actively operated satellites \cite{esasustainabilityreport}. Due to the high relative speeds involved, any collision would cause significant damage to assets. Satellites are therefore forced to perform \glspl{cam} when a critical close encounter is predicted.

Many satellites in crowded operational orbits primarily use low-thrust propulsion systems, such as ion engines, Hall-effect thrusters, and electrospray thrusters \cite{mazouffre_electric_2016}. These systems are more efficient than high-thrust propulsion, but they produce much weaker levels of thrust \cite{lev_technological_2019}; therefore, even small maneuvers such as \glspl{cam} (typically in the order of $1$-$10$~\si{cm/s}) can require extended thrusting windows. This limitation makes it important to start a satellite thrust arc adequately in advance, and can make optimizing these arcs a complex task due to the large parameter space \cite{de2022low}.

Existing literature on low-thrust \gls{cam} design primarily focuses on optimizing \gls{cam} thrust arcs to minimize fuel consumption, often resulting in early maneuvers \cite{hernando2021low}.
The \gls{cam} problem, therefore, is often framed as a \gls{fo} \gls{ocp}. Both analytical \cite{hernando2021low,Gonzalo2021,DeVittori2022} and numerical \cite{Patera2003, armellin2021collision, pavanello2026, pavanello2024recursive} approaches have been proposed in recent years to this end: a comprehensive review of them can be found in Ref. \cite{pavanello2025cammary}.
However, operational constraints often favor delaying \glspl{cam} as long as possible to leverage updated tracking data and reduce unnecessary maneuvers \cite{Schaus2025}. The reason for this discrepancy lies in the way the danger level of a satellite conjunction is quantified; typically in terms of the \gls{poc} or the deterministic \gls{md}, metrics computed at the \gls{tca} between the two objects. This uncertainty is regularly updated through ongoing tracking of the debris object: more recent tracking data yields lower positional uncertainty and, consequently, a more accurate assessment of conjunction danger and whether a \gls{cam} is needed. 

Therefore, following the research trail started in Ref. \cite{dell2024characterization}, this work shifts the focus from a \gls{fo} to a \gls{to} formulation: given a very short warning time, the interest is on the optimization of the latest initiation time for the \gls{cam}. Not only is this kind of optimization pertinent to collision alerts received very late, but it can also be used to determine whether there is enough time left to wait for the next batch of debris-tracking data to arrive and reassess whether a \gls{cam} is necessary.

The methodology is based on the idea of a \textit{backward-sweep}, wherein, starting from \gls{tca}, the low-thrust control profile is continuously updated by stacking actions backward until safety is reached. The control action is computed using a greedy approach similar to \cite{Patera2003} and \cite{pavanello2024recursive}, without the need to exactly solve the \gls{ocp}; thus, the method is named \gls{gto}. This can be done by expanding the dynamics' flow to higher order using \gls{da}: the collision metric at \gls{tca} becomes a polynomial function of the control actions applied at predefined time nodes. By substituting the greedy control computed during the previous step of the backward sweep, one can compute the greedy control of the current step and repeat the process until the safety threshold is reached. The use of \gls{da}, moreover, allows for the online update of the \gls{tca} shift due to the maneuver, as proven in previous works \cite{armellin2021collision,pavanello2024recursive}.

The remainder of this work is organized into four sections. In \cref{sec:prob}, the dynamics of the system are introduced, and the \gls{to} \gls{cam} problem is set. \cref{sec:back_sweep} explains the \gls{da}-based methodology to compute the approximated maps and establish the \gls{gto} solution method. Finally, in \cref{sec:results}, numerical results are shown, and conclusions are drawn in \cref{sec:conclusions}.

\section{CAM Problem Setting}
\label{sec:prob}
We consider the state of a controlled primary satellite at \gls{tca} described by a multivariate Gaussian random variable
\begin{equation}
\vec{X}_p(t_{CA}) = \begin{bmatrix} \vec{R}_p^\top \ \vec{V}_p^\top \end{bmatrix}^\top\sim \mathcal{N}_6(\vec{x},\vec{P}_p)
\end{equation}
where $\vec{R}_p\sim\mathcal{N}_3(\vec{r}_p,\vec{\Sigma}_p)$ and $\vec{V}_p\sim\mathcal{N}_3(\vec{v}_p,\vec{0}_3)$ are its position and velocity, respectively. All quantities here represent quantities at \gls{tca}.
The primary has a close encounter with a secondary (assumed uncontrolled) object
\begin{equation}
\vec{X}_s(t_{CA}) = \begin{bmatrix} \vec{R}_s^\top \ \vec{V}_s^\top \end{bmatrix}^\top\sim \mathcal{N}_6(\vec{x}_s,\vec{P}_s)
\end{equation}
where $\vec{R}_s\sim\mathcal{N}_3(\vec{r}_s,\vec{\Sigma}_s)$ and $\vec{V}_s\sim\mathcal{N}_3(\vec{v}_s,\vec{0}_3)$. For both the primary and secondary object, the velocity uncertainty is assumed to be null, in accordance with widely assessed \gls{cam} routines from the literature \cite{pavanello2025cammary}.
The dynamics of the mean state of the two objects are given by:
\begin{subequations}
\label{eq:StateDynamics}
\begin{align}
	\dot{\vec{x}}(t)   &= F(\vec{x}(t), t) + G(\vec{x}(t), t) \vec{u}(t), \label{eq:stateDyn_p} \\
	\dot{\vec{x}}_s(t) &= F(\vec{x}_s(t), t), \label{eq:stateDyn_s}
\end{align}
\end{subequations}
where $t\in [t_{alert}, \ t_{CA}]$; \( F(\cdot): \mathbb{R}^6 \times \mathbb{R} \rightarrow \mathbb{R}^6 \) denotes the drift, representing uncontrolled motion of the objects, with $t_{alert}$ and $t_{CA}$ the extremes (alert receive and \gls{tca}) of the time window of interest. \( F(\cdot) \) can include any dynamics, though in this work we limit the drift to two-body acceleration only, since it has been proven that perturbations have no significant effect on the short-term analysis of \glspl{cam} \cite{hernando2021low}. \( G(\cdot): \mathbb{R}^6   \times \mathbb{R} \rightarrow \mathbb{R}^{6 \times 3}\) denotes the influence of control on the dynamics, with \( \vec{u} \in \mathbb{R}^3 \) being the control, such that it is always at maximum throttle $||\vec{u}(t)||=u_{max}$, \( u_{max} \in \mathbb{R_+} \). In this work, we express $\vec{u}$ in \gls{rtn}-components, meaning that $G(\cdot)$ involves a frame transformation from \gls{rtn} to \gls{eci}, in which the state vector is expressed. It is worth noting that the propagation of the covariance is not of interest for the test case at hand, so no uncertainty propagation method is needed, and the two covariances are assumed as problem parameters at \gls{tca}. This is a valid assumption even under the hypothesis of variable \gls{tca}, as it has been shown in previous work that typical \gls{tca} shifts for minimum-time maneuvers are below $0.5$ \si{s}, a time frame in which the covariance is not expected to vary significantly \cite{pavanello2024recursive}.
The \gls{gto} method is only valid in the context of \textit{short-term encounters}, where the relative velocity between the two objects is high, and the conjunction event can be considered instantaneous at \gls{tca} \cite{patera2003general}. Therefore, the chosen danger metric is only assessed at \gls{tca}, rather than in an extended time interval. Since the \gls{gto} methodology tracks the \gls{tca} shift due to the maneuver, under the short-term approximation, it is ensured that the conjunction mitigation is effective.

\subsection{B-plane dynamics}

At \gls{tca}, one can define the multi-variate Gaussian random variable of the relative position as

\begin{equation}
\label{eq:rel_pos}
    \vec{R}(t_{CA}) = \vec{R}_p - \vec{R}_s \sim \mathcal{N}_3(\vec{r},\vec{\Sigma}),
\end{equation}
where $\vec{r} = \vec{r}_p - \vec{r}_s$ and $\vec{\Sigma} = \vec{\Sigma}_p + \vec{\Sigma}_s$. The relative velocity, which is deterministic, is defined as 
\begin{equation}
\label{eq:rel_vel}
    \vec{v}(t_{CA}) = \vec{v}_p - \vec{v}_s.
\end{equation}
Given the short-term nature of the encounter, it can be studied in the B-plane reference frame denoted as $\mathcal{B}$ \cite{pavanello2025cammary}. $\mathcal{B}$ is centered on the secondary object;
the $\eta$ axis aligns with the direction of the relative velocity of the primary with respect to the secondary, while the $\xi\zeta$ plane is perpendicular to the relative velocity axis \cite{armellin2021collision}.
Since the nominal \gls{tca} is the moment when the \gls{md} is minimum, by assumption, 
\begin{equation}
	\vec{r}(t_{CA})\cdot\vec{v}(t_{CA}) = 0,
\label{eq:TrackTCA}
\end{equation}
indicating that the relative position lies within the $\xi\zeta$ plane. From this point forward, the time argument of the variables will be shadowed, as all quantities are intended as referred to $t_{CA}$.
Therefore, for the computation of \gls{poc}, we will employ the projection of the relative state and its covariance on the $\xi\zeta$ plane: $\Vec{r}_\mathcal{B}\in\mathbb{R}^2$ and $\vec{\Sigma}_\mathcal{B}\in\mathbb{R}^{2\times2}$.
\begin{subequations}
\begin{align}
	& \Vec{r}_\mathcal{B} = \vec{Q} \ \vec{r},
	\label{eq:bPlaneProjr} \\
	& \vec{\Sigma}_\mathcal{B} = \vec{Q} \ \vec{\Sigma} \ \vec{Q}\transp.
	\label{eq:bPlaneProjP}
\end{align}
\label{eq:bPlaneProj}
\end{subequations}
where $\vec{Q}(\cdot):\mathbb{R}^3\times\mathbb{R}^3\rightarrow\mathbb{R}^{2\times 3}$ performs the rotation between the reference frame used for the propagation of the states (\gls{eci}) and the B-plane
\begin{equation}
\vec{Q} = \vec{Q}(\vec{v}_p,\vec{v}_s) = [\hat{u}_{\xi}\transp \ \hat{u}_{\zeta}\transp],
\label{eq:R}
\end{equation}
where $\hat{u}_{\xi}$ and $\hat{u}_{\zeta}$ are the unit vectors of the $\xi$ and $\zeta$ axes. The explicit expression of $\vec{R}$ can be found in reference \cite{armellin2021collision}.
The subsequent discussion assumes that all uncertainty is concentrated around the secondary object, while all mass is concentrated around the primary \cite{li_review_2022}.

A conservative way to compute \gls{poc} of the encounter is enveloping each of the bodies in a sphere with a radius equal to the largest dimension. The sum of the radii of the two spheres is the combined \gls{hbr}, which identifies the combined hard body circle $\mathbb{C}_\mathrm{HBR}$ on the B-plane.
To compute \gls{poc}, denoted as $P_C\in\mathbb{R}_\mathrm{+}$, the probability density function of the projection of the relative position onto the B-plane is integrated over the combined hard-body circle, as follows
\begin{equation}
\begin{split}
\pc = & \frac{1}{(2 \pi)^{3/2}\sqrt{\mathrm{det}\big(\boldsymbol{P}_\mathcal{B}\big)}} \cdot \\ & \cdot \iint_{\mathbb{C}_\mathrm{HBR}} \mathrm{exp}\left(-\frac{\smd}{2}\right) \mathrm{d} A,
\end{split}
\label{eq:pceq}
\end{equation} 
where $\smd\in\mathbb{R_+}$ is the \gls{smd} of the relative position in the B-plane and is defined as 
\begin{equation}
    \label{eq:smd}
    \smd = \vec{r}_\mathcal{B}\transp\vec{\Sigma}_\mathcal{B}^{-1}\vec{r}_\mathcal{B}.
\end{equation}
\gls{smd} provides a normalized measure of separation that incorporates both the relative distance and the positional uncertainties.

Multiple approaches have been proposed to approximate the integral in \cref{eq:pceq} \cite{li_review_2022}. In general, for a given \gls{hbr}, \gls{poc} is expressed as a function of the relative position, and the combined covariance in the B-plane at \gls{tca}
\begin{equation}
    \pc = \pc(\vec{r}_\mathcal{B},\vec{\Sigma}_\mathcal{B}).
    \label{eq:generalPoc}
\end{equation}
However, in this work, we make use of the analytical expression formulated by Chan \cite{Chan2008Spacecraft}, which allows for the inversion of the formula to obtain exact isoprobability curves \cite{de_vittori_combined_2026}.
Therefore, without loss of accuracy, it is possible to convert the safety condition from \gls{poc} to \gls{smd}, which is a more convenient metric to represent in \gls{da} since it is quadratic in the relative position.

An alternative approach to evaluating the safety of a close encounter is to disregard information about covariance, aiming for a deterministic metric. This approach is popular when the estimation of the uncertainty is unreliable \cite{symonds2014operational}. In this case, the used metric is the \gls{md}, defined, at \gls{tca}, as 
\begin{equation}\label{eq:EucDis}
	d_{miss} = \| \vec{r}_\mathcal{B} \|.
\end{equation}
From now on, for simplicity of exposition, $d_{(\cdot)}$ represents the danger metric of a conjunction, without specifying the exact metric choice (\gls{md} or \gls{smd}).

\subsection{Optimal Control Problem}

Regardless of the danger metric used (\gls{md} or \gls{smd}), the goal of a \gls{cam} is to bring it over a safe threshold at \gls{tca}, indicated with $\gamma\in\mathbb{R_+}$. The objective investigated in this paper is to maximize the time available before initiating the \gls{cam}, which is equivalent to minimizing $-t_0$, where $t_0\in[t_{alert}, \ t_{CA}]$ is the start time of the maneuver (assumed positive):
\begin{subequations}\label{eq:OriginalProblemDA}
	\begin{align}
		\min_{\vec{u}} \quad & -t_{0} \label{eq:ocp_obj} \\ 
		\text{s. t.} \quad & \dot{\vec{x}} = F(\vec{x},t) + G(\vec{x},t)\vec{u} \label{eq:ocp_dyn} \\ 
		\quad & \dot{\vec{x}}_s = F(\vec{x}_s,t) \label{eq:ocp_dyn_s} \\ 
		\quad & \vec{x}(t_0) = \mathcal{F}_{CA}(\vec{x}(\bar{t}_{CA}),t_{CA},t_0) \label{eq:ocp_x0} \\
        \quad & t_{CA} = \bar{t}_{CA} + f(\vec{x}(t_0), \vec{x}_s(\bar{t}_{CA}), \vec{u}(t)) \label{eq:ocp_tca} \\ 
		\quad & d_{(\cdot)}(t_{CA}) \geq \gamma \label{eq:ocp_metric} \\
        \quad & ||\vec{u}|| = u_{max}, \label{eq:ocp_umax}
        \end{align}
\end{subequations}
where the time arguments of states and control have been neglected, where needed for clarity. In Problem~\eqref{eq:OriginalProblemDA}, \cref{eq:ocp_obj} is the \gls{to} objective function; \cref{eq:ocp_dyn,eq:ocp_dyn_s} are the dynamics constraints of the primary and secondary spacecraft, respectively, from \cref{eq:StateDynamics}; in \cref{eq:ocp_x0}, the initial state \( \vec{x}(t_0) \) is determined through the backward propagation of the state of the primary satellite from the nominal \gls{tca}, assumed known, e.g., from a \gls{cdm}:
\begin{equation}
\label{eq:flow}
    \mathcal{F}_{CA} = \int_{t_{CA}}^{t_0} F \left(\vec{x}(t)\right) dt
\end{equation} 
where $\mathcal{F}_{CA}(\cdot)$ is the backward-integrated flow of $F(\cdot)$; \cref{eq:ocp_tca} accounts for the change in \gls{tca} due to the maneuver, where ${f(\cdot):\mathbb{R}^6\times\mathbb{R}^6\times\mathbb{R}^3\rightarrow\mathbb{R}}$ is the generic function that identify the closest approach condition, $t_{CA}$ is the shifted value and $\bar{t}_{CA}$ is the nominal value; \cref{eq:ocp_metric} is the \gls{smd} or \gls{md} constraint, computed at the shifted \gls{tca}. Finally, \cref{eq:ocp_umax} is the constraint on the thrust to always be at maximum throttle $u_{max}\in\mathbb{R_+}$.

We stress the importance of accounting for the variation of \gls{tca} with the control input for accuracy, since the close approach can move due to the maneuver by up to a second, in certain conditions \cite{pavanello2024recursive}. At the point of geometric closest approach, the relative position and velocity vectors are orthogonal (see \cref{eq:TrackTCA}).
In typical short-term encounters, the relative distance between the two objects changes more than the relative velocity as a result of the \gls{cam}. In general, \cref{eq:TrackTCA} will no longer be valid for the nominal value of \gls{tca} without any control, and the time at which the relative distance is minimal changes as a result of applied control. While this orthogonality condition holds strictly only when considering the \gls{md} as a danger metric, we use it as a first-order approximation for the critical time of any given danger metric, as the geometric closest approach always aligns closely with the time of maximum danger when looking at \gls{poc} or \gls{smd}. 

\section{Backward Sweep Approach}
\label{sec:back_sweep}

To solve the optimization problem in Problem \eqref{eq:OriginalProblemDA}, we employ a numerical approach based on a backward sweep in time, enabled by \gls{da}. Instead of solving the continuous optimization problem directly, we discretize the time interval and iteratively determine the optimal control action that maximizes the rate of change of the chosen danger metric. We choose to start this sweep from $t_{CA}$ rather than $t_0$ for two reasons. (i) The effect of control in a certain time interval on the danger metric at \gls{tca} is dependent on the control in future time intervals due to the presence of natural drift; when performing a backward rather than a forward sweep, these effects are already known. (ii) Since the optimization problem concerns finding the last possible moment to initiate thrust, it is sensible to perform a backward sweep starting from \gls{tca} and assess continuously whether the safety criterion has already been reached.

In the context of the optimization problem shown in Problem~\eqref{eq:OriginalProblemDA}, we assume that a \gls{cdm} is received at $t_{alert}$, such that $t_{alert} < t_0 < t_{CA}$. It is not strictly known beforehand whether the time of receiving the collision alert is before the last possible moment of initiating thrust, though typically $t_{alert}$ is in the order of days before \gls{tca}.
The continuous time window ${t\in\mathbb{R}_{[t_{alert}, \ t_{CA}]}}$ is discretized into a time grid with $N+1$ nodes ${t_i\in\{t_{alert}, t_1, ..., \bar{t}_{CA}\}}$. 
Starting with the last time interval $[t_{N - 1}, \bar{t}_{CA}]$ in the backward sweep, the control acceleration is determined that maximizes the rate of change of the relevant danger metric. Note that this does not exactly solve the problem posed in Problem~\eqref{eq:OriginalProblemDA}, as gradient maximization of the metric is not equivalent to maximizing $t_0$ because it does not account for higher order terms.

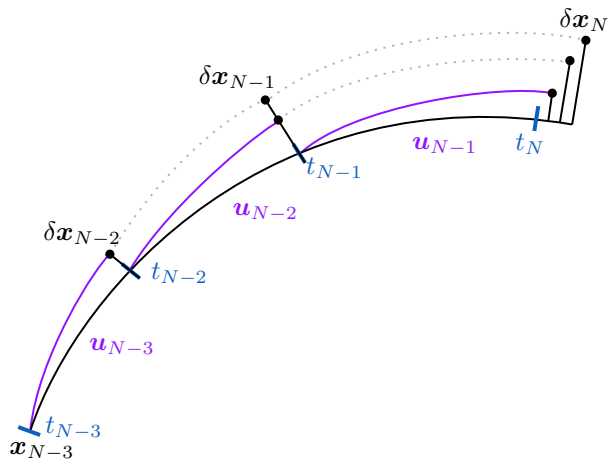
\begin{figure}[b!]
    \centering
    \input{Figures/back_sweep.tex}
	\caption{Illustration of the backward sweep procedure. The nominal trajectory is shown in black with the adjusted trajectory in purple.}
	\label{fig:DiagramOfGeometry}
\end{figure}

Tracking $\vec{r}_{\mathcal{B}}$ as a function of the applied control in the time interval $[t_{N - 1}, \bar{t}_{CA}]$, the change in \gls{md} at \gls{tca} (and through that, the change in \gls{smd} and \gls{poc}) is tracked. It can then be assessed whether the resulting change is sufficient to overcome the safety threshold for the relevant danger metric. If not, the procedure is repeated for the time interval $[t_{N-2}, t_{N-1}]$, continuing the backward sweep until the safety threshold $\gamma$ is reached. If at some point the safety threshold is reached, it can be concluded that $t_0$ lies in the time interval analyzed last in the backward sweep. 
The idea of the \gls{gto} method is sketched out in \autoref{fig:DiagramOfGeometry}, with a nominal trajectory in blue along discretized time steps $t_{j}$, $j \in \{0, \ ... \ , 3\}$

In the next discussion, it is described in mathematical detail how exactly the control is determined and how the effects of state and control perturbations are coupled with the danger metric and \gls{tca} behavior. 

\subsection{Differential Algebra}
The numerical propagation of the spacecraft's state and its sensitivities is accomplished using \gls{da}. In the context of astrodynamics and this paper, \gls{da} refers to a powerful computational technique for automatic differentiation that computes high-order Taylor series expansions of functions \cite{armellin2010asteroid, morselli2014high}. The core idea is to replace standard floating-point arithmetic with operations on multivariate Taylor polynomials. Each variable (e.g., a state component) is represented not as a single number, but as a polynomial that encodes its value and its partial derivatives with respect to a set of independent monomials, representing initial perturbations.

Let $\vec{y} \in \mathbb{R}^m$ be a vector of independent perturbed variables, such as initial state deviations and control actions, expanded around a nominal point $\bar{\vec{y}}$. Any function $g(\vec{y})$ is represented in the \gls{da} framework as its $k^\text{th}$-order Taylor expansion:
\begin{equation}
	g(\vec{y}) \approx \mathcal{T}^k_{g(\vec{y})}(\vec{y}) = g(\bar{\vec{y}}) + \sum_{k=1}^{n} \frac{1}{k!}\frac{\partial^k g}{\partial \vec{y}^k}(\bar{\vec{y}})(\vec{y}-\bar{\vec{y}})^k,
\end{equation}
where $\bar{\vec{y}} \in \mathbb{R}^m$ is the constant part of the expanded variable (expansion point).
When an operation, be it any arbitrary function, is required on these \gls{da} variables, the chain rule is automatically applied using the Taylor expansion of transcendental functions and algebraic operations between polynomials to deliver a composed Taylor expansion.

In the \gls{gto} methodology, \gls{da} is used to construct a $k^{\text{th}}$-order polynomial map (typically $2^{\text{nd}}$-order), $\mathcal{T}^k$, representing the flow of the state dynamics described earlier. Working in the discretized time grid with nodes $t_i$, this map expresses the deviation of the state at a future time, $\delta\vec{x}(t_{i+1})$, as an explicit polynomial function of the state perturbation at a previous time, $\delta\vec{x}(t_i)$, and the control perturbation $\delta\vec{u}(t_i)$ applied over the interval $[t_i, t_{i+1}]$:
\begin{equation}
    \label{eq:daDyn}
	\delta\vec{x}(t_{i+1}) = \mathcal{T}^k_{\delta\vec{x}(t_{i+1})}(\delta\vec{x}(t_i), \delta\vec{u}(t_i)).
\end{equation}
This approach allows for the efficient and highly accurate propagation of not only the state, but also its complex, nonlinear sensitivities to perturbations, which is fundamental to the backward sweep algorithm.

\subsection{Backward Sweep of the First Interval}
The nominal state at $t_{N-1}$ is obtained by back-propagating the nominal value of $\vec{x}_{N}$ from the nominal \gls{tca}
\begin{equation}
    \label{eq:backProp}
    \bar{\vec{x}}_{N-1} = \mathcal{F}_{N}(\bar{\vec{x}}_{N}, \bar{t}_{CA},t_{N-1}),
\end{equation}
where $\mathcal{F}_{N}(\cdot)$ indicates the integrated dynamics flow from $t_N = \bar{t}_{CA}$ to $t_{N-1}$.
The state of the primary satellite at time \( t_{N-1} \) is now expanded around its nominal value
\begin{equation}
	\vec{x}(t_{N-1}) = \bar{\vec{x}}(t_{N-1}) + \delta \vec{x}(t_{N-1}),
\end{equation}
so that the state deviation at node $N$ is a Taylor polynomial in the state and control deviation of the previous node. That is to say, $\vec{x}(t_{N-1})$ can be described by the sum of a nominal state (\( \bar{\vec{x}}(t_{N-1}) \)) and \gls{da} perturbations (\( \delta \vec{x}(t_{N-1}) \)) caused by the control \( \delta \vec{u} \) of previous time steps.
The control is also expanded around the nominal null value ($\bar{\vec{u}}_{N-1} = 0$), so that $\vec{u}_{N-1} = \delta \vec{u}_{N-1}$, constant over the interval $[t_{N-1}, t_N]$. 
Lastly, \gls{tca} is also expanded around its nominal value $\bar{t}_{CA}$, so that the effect of the maneuver on it can be tracked:
\begin{equation}
    t_{CA} = \bar{t}_{CA}  + \delta t_{CA}
\end{equation}

Following \cref{eq:StateDynamics} and using a high-fidelity integration scheme like Runge-Kutta 7-8, the following discretized dynamics are obtained, coupling the state of the primary node $N$ to the state at node $N-1$:
\begin{equation}  
\begin{aligned}
	\vec{x}_{N} & = \mathcal{F}_{N-1}(\vec{x}_N, t_{N-1},  t_N) \ + \\ &+ \mathcal{G}_{N-1}(\vec{x}_{N-1}, t_{N-1},  t_N) \vec{u}_{N-1}.  \label{eq:StateDynamicsDiscretised_p}
	\end{aligned}
\end{equation}
The same can be done for the secondary object
\begin{equation}
    \vec{x}_{s,N} = \mathcal{F}_{N-1}(\vec{x}_{s,N-1}, t_{N-1},  t_N),
    \label{eq:StateDynamicsDiscretised_s}
\end{equation}
where $\vec{x}_i = \vec{x}(t_i)$, and $\mathcal{F}_{i}(\cdot)$ and $\mathcal{G}_{i}(\cdot)$ represent the flow of the drift and controlled dynamics from time $t_i$ to $t_{i+1}$, computed as in \cref{eq:flow}. 
The control $\vec{u}_i = \vec{u}(t_i)$ is expressed in \gls{rtn} components within each time interval. Note that in this time interval, the length is not fixed but rather dependent on $\delta t_{CA}$; without control, $\delta t_{CA}$ is zero, but under the influence of control in each time interval, $\delta t_{CA}$ will change and accumulate throughout the backward sweep. Since the backward-sweep does not update with a forward propagation when successive time intervals are included, the control applied to the last time interval only acts on the updated \gls{tca} computed during the first sweep.

Since the variables in \cref{eq:StateDynamicsDiscretised_p,eq:StateDynamicsDiscretised_s} are \gls{da} variables, the equations can be equivalently written in the form of polynomials, such as in \cref{eq:daDyn}
\begin{subequations}
\label{eq:StateTaylors}
\begin{align}
	\vec{x}_{N} &= \mathcal{T}_{\vec{x}_N}(\delta \vec{x}_{N-1}, \delta \vec{u}_{N-1}, \delta t_{CA})  \\
	\vec{x}_{s,N} &= \mathcal{T}_{\vec{x}_{s,N}}(\delta t_{CA}),
\end{align}
\end{subequations}
where the state of the secondary object is unaffected by control, and the Taylor polynomials describing its state in \cref{eq:Taylors} are therefore only dependent on $\delta t_{CA}$. Splitting the states into positions and velocities, one can write four multivariate Taylor polynomials
\begin{subequations}
\begin{align}
	& \vec{r}_{p,CA} = \mathcal{T}_{\vec{r}_{p,CA}}(\delta \vec{x}_{N-1}, \delta \vec{u}_{N-1}, \delta t_{CA}) \label{eq:x1} \\
    & \vec{v}_{p,CA} = \mathcal{T}_{\vec{v}_{p,CA}}(\delta \vec{x}_{N-1}, \delta \vec{u}_{N-1}, \delta t_{CA}) \label{eq:x2} \\
    & \vec{r}_{s,CA} = \mathcal{T}_{\vec{r}_{s,CA}}(\delta t_{CA}) \label{eq:x3} \\
    & \vec{v}_{s,CA} = \mathcal{T}_{\vec{v}_{s,CA}}(\delta t_{CA}).\label{eq:x4} 
\end{align}
        \label{eq:Taylors}
\end{subequations}
Consequently, the relative position and velocity at \gls{tca} can also be expressed as Taylor polynomials, by subtracting \cref{eq:x3} from \cref{eq:x1} and \cref{eq:x4} from \cref{eq:x2}.
\begin{subequations}
\label{eq:relTaylor}
    \begin{align}
        & \vec{r}_{CA} = \mathcal{T}_{\vec{r}_{CA}}(\delta \vec{x}_{N-1}, \delta \vec{u}_{N-1}, \delta t_{CA}) \label{eq:relPosTaylor} \\
        & \vec{v}_{CA} = \mathcal{T}_{\vec{v}_{CA}}(\delta \vec{x}_{N-1}, \delta \vec{u}_{N-1}, \delta t_{CA}),\label{eq:relVelTaylor}
    \end{align}
\end{subequations}
where $\vec{r}_{CA} = \vec{r}(t_{CA})$ and $\vec{v}_{CA} = \vec{v}(t_{CA})$. Lastly, the rotation matrix from \gls{eci} to \gls{rtn} in \cref{eq:R} is also expanded around the state and time deviation at \gls{tca}
\begin{equation}
\vec{Q} = \mathcal{T}_{\vec{Q}}(\delta\vec{v}_p, \delta t_{CA}) 
\label{eq:TaylorR}
\end{equation}
where $\delta\vec{v}_p$ is directly derived from the control perturbation history and $\delta t_{CA}$ is the only \gls{da} perturbation affecting the velocity of the secondary at \gls{tca}.
The perturbation $\delta t_{CA}$ represents the \textit{total} shift in \gls{tca}, resulting from all control actions coming from the initiation of the maneuver up to conjunction. Here, the fact that low-thrust propulsion is used and consequently that the state deviation remains small \cite{pavanello2025cammary}, enables us to perform the expansion from the nominal state instead of the displaced state. Displacement in the tangential direction could still be significant even with low-thrust propulsion, but the limited time window considered and the bounds given by the safety threshold $\gamma$ prevent inaccuracies in this regard. This is important because it allows us to compute the Taylor expansions that map the states from $t_i$ to $t_{i+1}$ all at once: in the backward sweep, the dynamics are not re-propagated when a new segment is included, resulting in more computationally efficient solutions.


The dependency on $\delta t_{CA}$ in Eqs.~\eqref{eq:Taylors} and ~\eqref{eq:TaylorR} can be dropped by means of a polynomial partial inversion technique. Let us consider the close approach orthogonality condition from \cref{eq:TrackTCA}: according to \cref{eq:relTaylor}, its Taylor expansion reads
\begin{equation}
\mathcal{T}_{\vec{r}_{CA} \cdot \vec{v}_{CA}}(\delta\vec{x}_{N-1},\delta\vec{u}_{N-1},\delta t_{CA}) = 0.
\label{eq:paramImplicit}
\end{equation}
The parametric implicit constraint in \cref{eq:paramImplicit} can be solved for $\delta t_{CA}$ using the method proposed in Ref. \cite{armellin2010asteroid}, which is also used in Refs. \cite{armellin2021collision,pavanello2024recursive} in the context of \gls{cam} optimization
\begin{equation}
\begin{aligned}
\delta t_{CA} &= \mathcal{T}_{\delta t_{CA}}(\vec{r}_{CA} \cdot \vec{v}_{CA} = 0,\delta\vec{x}_{N-1},\delta\vec{u}_{N-1}) \\
&= \mathcal{T}_{\delta t_{CA}}(\delta\vec{x}_{N-1},\delta\vec{u}_{N-1}),
\label{eq:deltaT}
\end{aligned}
\end{equation}
Modern \gls{da} libraries, such as DACE\footnote{\url{https://github.com/dacelib/dace}}, provide built-in functions to perform this inversion automatically.

Substituting the polynomial from \cref{eq:deltaT} into \cref{eq:relTaylor,eq:TaylorR}, one can eliminate the explicit dependency of the conjunction variables on the \gls{tca} perturbation, yielding ${\vec{r}_{CA} = \mathcal{T}_{\vec{r}_{CA}}(\delta\vec{x}_{N-1},\delta\vec{u}_{N-1})}$ and ${\vec{Q} = \mathcal{T}_{\vec{Q}}(\delta\vec{x}_{N-1},\delta\vec{u}_{N-1})}$. 
These are substituted into \cref{eq:bPlaneProjr,eq:bPlaneProjP}, so the relative position and positional covariance in the B-plane become functions of the state deviation and the control at $t_{N-1}$  alone:
\begin{subequations}
\label{eq:r_and_sigma}
\begin{align}
	& \vec{r}_\mathcal{B} =  \mathcal{T}_{\vec{r}_\mathcal{B}}(\delta\vec{x}_{N-1},\delta\vec{u}_{N-1}),
\label{eq:expBplane} \\
    & \vec{\Sigma}_\mathcal{B} =  \mathcal{T}_{\vec{\Sigma}_\mathcal{B}}(\delta\vec{x}_{N-1},\delta\vec{u}_{N-1}). \label{eq:covBplane}
\end{align}
\end{subequations}
Consequently, the danger metric at $t_{CA}$  can be approximated by a Taylor polynomial according to either \cref{eq:smd} or \cref{eq:EucDis}, since they are only dependent via the chain rule on the expressions in \cref{eq:r_and_sigma}
\begin{equation}\label{eq:DMTaylor}
	d_{(\cdot)} = \mathcal{T}_{d_{(\cdot)}}(\delta \vec{x}_{N-1}, \delta \vec{u}_{N-1})
\end{equation}
\subsection{Greedy control optimization}

To determine the control perturbation $\delta \vec{u}$ maximizing $t_0$, we investigate the first order term of the danger metric $d_{(\cdot)}$:
\begin{equation}\label{eq:DMExpansion}
	d_{(\cdot)} = \bar{d}_{(\cdot)} + \frac{\partial d_{(\cdot)}}{\partial \vec{u}_{N-1}} \delta \vec{u}_{N-1}
\end{equation}

In low-thrust propulsion scenarios, the limited thrust magnitude ($u_{max} \ll 1$) necessitates precise targeting of the rate of change of the danger metric to achieve optimal collision avoidance. To maximize the linear rate of change of the danger metric, the control acceleration is chosen to be aligned with the gradient, as also done by Patera and Peterson \cite{patera2003general}, and by the authors in Refs. \cite{pavanello2024recursive,Pavanello2025efficient}:
\begin{equation}\label{eq:Control}
	\bar{\vec{u}}_{N-1} = u_{max} \frac{ \frac{\partial d_{(\cdot)}}{\partial \vec{u}_{N-1}}}{|| \frac{\partial d_{(\cdot)}}{\partial \vec{u}_{N-1}}||},
\end{equation}
where the overline indicates that the control is computed and not a \gls{da} variable.
Note that \cref{eq:Control} does not exactly solve the optimization problem, and results in a sub-optimal maneuver, as it aligns the thrust with the local gradient of the metric for computational efficiency. 

The result from \cref{eq:Control} is plugged into the polynomial map \cref{eq:DMTaylor}, so that the metric variation is only a function of the state deviation at node $N-1$
\begin{equation}
\begin{aligned}  
    d_{(\cdot)}^{(1)} &= \mathcal{T}_{d_{(\cdot)}}(\delta \vec{x}_{N-1},\delta\vec{u}_{N-1} = \bar{\vec{u}}_{N-1}) \\ &= \mathcal{T}_{d_{(\cdot)}^{(1)}}(\delta \vec{x}_{N-1})  \label{eq:d_plugged_1} 
    \end{aligned}
\end{equation}
where the apex $(1)$ indicates that this multivariate polynomial is obtained by evaluating the previous one in the constant control action. It should be stressed that no new \gls{da} maps have been built via function evaluations, but the new maps have simply been obtained by evaluating the original map using the constant control action, enhancing computational efficiency. The first-order approximation of the effect of the control on the metric, therefore, is the constant part of the polynomial in \cref{eq:d_plugged_1}: $\bar{d}_{(\cdot)}^{(1)}$.

\subsection{Backward Sweep at Previous Time Steps}

If $\bar{d}_{(\cdot)}^{(1)}$ has not yet reached the safety threshold $\gamma$, the backward sweep must continue. To this end, $\delta \vec{x}_{{N-1}}$ must be formulated as a Taylor polynomial of perturbations in the time interval $[t_{N-2},t_{N-1}]$, and so forth, cascading the effect of the perturbations up to the necessary time $t_{N-i}$. 

So, let us now consider the general case in which the maps have been progressively evaluated through cascaded effects up to node $N-i+1$. The Taylor map of the metric in $\delta\vec{x}_{N-i}$, analogous to \cref{eq:d_plugged_1}, reads
\begin{equation}
\begin{aligned}  
    d_{(\cdot)}^{(i-1)} &= \mathcal{T}_{d_{(\cdot)}^{(i-2)}}(\delta \vec{x}_{N-i+1},\delta\vec{u}_{N-i+1} = \bar{\vec{u}}_{N-i+1}) \\ &= \mathcal{T}_{d_{(\cdot)}^{(i-1)}}(\delta \vec{x}_{N-i+1}) \label{eq:d_plugged_i} 
    \end{aligned}
\end{equation} 

The back-propagated nominal state at $t_{N-i}$ is computed as
\begin{equation}
    \label{eq:backPropi}
    \bar{\vec{x}}_{N-i} = \mathcal{F}_{N-i+1}(\bar{\vec{x}}_{N-i+1}, t_{N-i+1}, t_{N-i}),
\end{equation}
and the state is expanded around the nominal value 
\begin{equation}
\label{eq:exp_state}
    \vec{x}_{N-i} = \bar{\vec{x}}_{N-i} + \delta \vec{x}_{N-i}.
\end{equation}
Therefore, the state deviation at node $N-i+1$ is written as a Taylor polynomial in the state and control deviation of the previous node, expanding \cref{eq:StateDynamicsDiscretised_p}
\begin{equation}
    \label{eq:TaylorsBack_i}
	\delta\vec{x}_{N-i+1} = \mathcal{T}_{\delta\vec{x}_{N-i+1}}(\delta \vec{x}_{N-i}, \delta \vec{u}_{N-i}).
\end{equation}
Importantly, from the second section of the sweep onward, the propagation time intervals are fixed and not variable as in the first sweep.

Substituting \cref{eq:TaylorsBack_i} inside the argument of the right side of \cref{eq:d_plugged_i}, one gets
\begin{equation}
\label{eq:ddot}
    d_{(\cdot)}^{(i\mathrm{-})}  = \mathcal{T}_{d_{(\cdot)}^{(i\mathrm{-})}} (\delta \vec{x}_{N-i}, \delta \vec{u}_{N-i}),
\end{equation} 
where the contribution of the controls computed in the previous steps is included in the constant part of the polynomials; the apex $(i\mathrm{-})$ indicates that this map is before the control evaluation, and is analogous to the original of \cref{eq:DMTaylor}.
For the generic node $N-i$, the greedy control is computed as in \cref{eq:Control}, by aligning it with the gradient of the metric
\begin{equation}\label{eq:Control_i}
	\delta \vec{u}_{N-i} = u_{max} \frac{ \frac{\partial d_{(\cdot)}^{(i\mathrm{-})}}{\partial \vec{u}_{N-i}}}{|| \frac{\partial d_{(\cdot)}^{(i\mathrm{-})}}{\partial \vec{u}_{N-i}}||}.
\end{equation}
At this point, the map is evaluated again using the computed control, yielding a pure dependency on the state deviation
\begin{equation}
    \label{eq:dend}
    d_{(\cdot)}^{(i)} = \mathcal{T}_{d_{(\cdot)}^{(i)}}(\delta \vec{x}_{N-i})   
\end{equation}

The constant part $\bar{d}_{(\cdot)}^{(i)}$ is checked, and the procedure is repeated for the previous time step if necessary, from \cref{eq:backPropi} to \eqref{eq:dend}, until the constant part of the metric exceeds the threshold. Moreover, the same passages that are taken to evaluate the maps of the metric can be used to progressively update the map of other relevant quantities, such as, most interestingly, the \gls{tca} change coming from \cref{eq:deltaT} or the relative position in the B-plane from \cref{eq:expBplane}. Note that these maps should also be evaluated using the controls computed via \cref{eq:Control_i}, not using controls in the maximum direction of change for themselves.

\subsection{Refinement of the last time interval}
Once the backward sweep reaches a value of the metric that is over the threshold, it is possible to apply a first-order expansion on the duration of the last time interval to approximately reach the exact threshold. Let us suppose that at the $j^\text{th}$ sweep the condition is satisfied: $\bar{d}_{(\cdot)}^{(j)} > \gamma$. A new \gls{da} variable is introduced, $\delta \tau\in\mathbb{R_+}$, which acts on the propagation time: the nominal state at $t_{N-j+1}$ is back-propagated taking into account this time perturbation:
\begin{equation}
    \label{eq:backProp_j}
    \vec{x}_{N-j} = \mathcal{F}_{N-j+1}(\bar{\vec{x}}_{N-j+1}, t_{N-j+1}, t_{N-j} + \delta \tau),
\end{equation} 
$\vec{x}_{N-j}$ is now a Taylor polynomial in $\delta \tau$ only. A new forward propagated state is obtained using the constant control found in the last sweep using \cref{eq:Control_i}
\begin{equation}
    \label{eq:forwProp_j}
    \begin{aligned}
    \vec{x}_{N-j+1} &  = \mathcal{F}_{N-j}(\vec{x}_{N-j}, t_{N-j}  + \delta \tau, t_{N-j+1}) + \\ & + \mathcal{G}_{N-j}(\vec{x}_{N-j}, t_{N-j} +\delta \tau, t_{N-j+1}) \vec{u}_{N-j}.
    \end{aligned}
\end{equation}
Therefore, $\vec{x}_{N-j+1}$ can be written as a polynomial in $\delta \tau$
\begin{equation}
    \label{eq:T_j}
    \vec{x}_{N-j+1} = \mathcal{T}_{\vec{x}_{N-j+1}}(\delta \tau)
\end{equation}
Using this and through the cascaded effect explained in the previous section, one can now obtain the first-order Taylor polynomial of the metric as a function of the time interval reduction
\begin{equation}
    d_{(\cdot)}^{(j)} = \mathcal{T}^1_{d_{(\cdot)}^{(j)}}(\delta\tau) = \bar{d}_{(\cdot)}^{(j)} + \sigma \frac{\partial d_{(\cdot)}^{(j)}}{\partial\tau}
\end{equation}
where $\sigma\in\mathbb{R_-}$ is the first-order coefficient of the polynomial. $\sigma$ should be negative since a reduction in the thrusting time (i.e., an increase in $\tau$) should result in a reduction of the metric.

In a first-order approximation, one can now find the required $\tau\in\mathbb{R_+}$ to subtract from the time interval $[t_{j-1}, \ t_j]$ to correctly reach the metric threshold:
\begin{equation}
    \label{eq:tau}
    \tau = \frac{\gamma - \bar{d}_{(\cdot)}^{(j)}}{\sigma}
\end{equation}
Evaluating $d_{(\cdot)}^{(j)}$ in  $\tau$, now allows to find the adjusted metric. Analogously, the \gls{tca} shift or the relative position in the B-plane can be written as polynomials in $\delta \tau$, and evaluating them in $\tau$ gives the adjusted value for a lower thrust duration.

\section{Results}
\label{sec:results}
The \gls{gto} algorithm is tested on an ESA’s Collision Avoidance Challenge, available for download, which includes the \glspl{cdm} of 2170 \gls{leo} conjunctions.\footnote{\url{https://github.com/arma1978/conjunction}} The dataset includes synthetic conjunctions in which the primary satellite is on a quasi-circular orbit, which is typical for \gls{leo} spacecraft.
Both results using \gls{smd} and \gls{md} are shown, in terms of optimality (maneuver initiation time), accuracy (validation errors), and computational efficiency (runtime). 
The \gls{gto} method is then compared with the \gls{ocp} formulation from Ref. \cite{dell2024characterization} to assess the validity of the sub-optimal solution obtained via the greedy control.
Lastly, it is compared with the solution of a \gls{fo} problem to assess the sub-optimality of the formulation in terms of fuel consumption.
The algorithm is run on an AMD Ryzen 9 6900HS at 3300 \si{Mhz} using MATLAB R2025b.
In \cref{tab:params}, the simulation parameters, valid throughout the 2170 cases, are shown.\footnote{The code and datasets developed for this work are available for download at \url{https://github.com/FrankdeVeld/CAM_DA}}

\begin{figure}[tb!]
    \centering
    \includegraphics[width=0.9\linewidth]{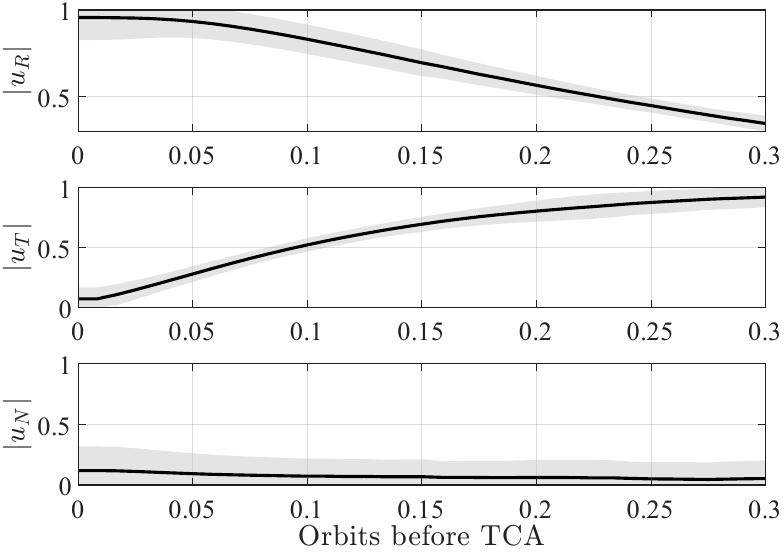}
    \caption{Statistics of the control profile for the \gls{smd} case.}
    \label{fig:control_smd}
\end{figure}

\begin{figure}[tb!]
    \centering
    \includegraphics[width=0.9\linewidth]{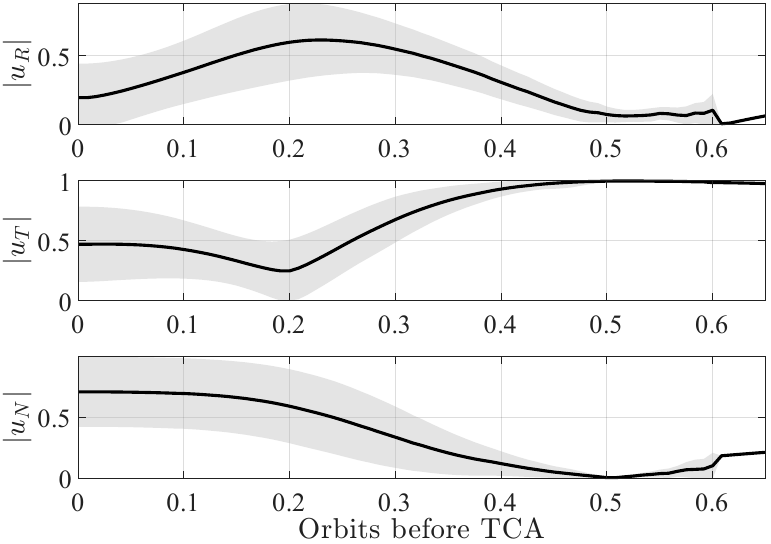}
    \caption{Statistics of the control profile for the \gls{md} case.}
    \label{fig:control_md}
\end{figure}

In \cref{fig:control_smd,fig:control_md}, the computed control evolution is shown for the \gls{smd} and the \gls{md} cases, respectively: a Gaussian fit is performed node-wise using the 2170 available samples from the dataset.
In these plots, the black lines indicate the mean values and the gray areas the 1-$\sigma$ deviation computed over the absolute value of the realizations;\footnote{The 1-$\sigma$ deviations are cropped at $0$ and $1$ not to represent unrealistic values.} this allows us to avoid an averaging of the mean values around $0$. It is clear that, as expected from discussions in the literature \cite{hernando2021low,pavanello2025cammary}, to decrease \gls{poc}, a radial thrust is preferable when close to the conjunction, while the greedy thrust increasingly moves to the tangential direction until it becomes purely tangential, on average, around $0.3$ orbits before \gls{tca}; the normal direction remains on average unused. 
The typical thrust to increase the \gls{md} shown in \cref{fig:control_md} follows a more unexpected profile, with a maximum value of radial thrust around $0.2$ orbits before \gls{tca} and a non-negligible out-of-plane component up to $0.4$ orbits before \gls{tca}; also in this case, the tangential direction becomes predominant when further from the conjunction. These behaviors are consistent with the evolution of the norm of the gradient of the metric, which is shown in Fig. 7 of Ref. \cite{pavanello2025cammary}.

\begin{table}[tb!]
\centering
    \caption{Simulation parameters.}
    \label{tab:params}
    \begin{tabular}{cccccc}
        \hline
        $u_{max}$ & $\gamma_{md}$ & $\gamma_{PoC}$ & $t_{alert}$ & \gls{da} order & N/orbit  \\  \hline
        $0.375$~\si{mm/s^2}                 & $2$~\si{km}                    & $10^{-6}$      & 1 orbit     & $2$ & $120$ \\
        \hline
    \end{tabular}
\end{table}

\begin{figure}
    \centering
    \includegraphics[width=0.9\linewidth]{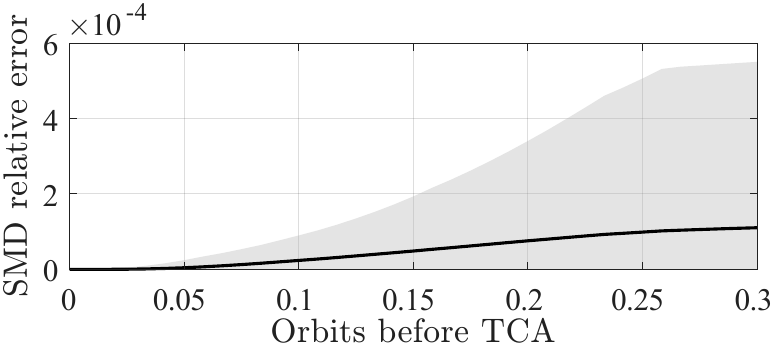}
    \caption{Statistics of the evolution of the relative \gls{smd} error against a forward-propagation validation.}
    \label{fig:error_smd}
\end{figure}

\begin{figure}
    \centering
    \includegraphics[width=0.9\linewidth]{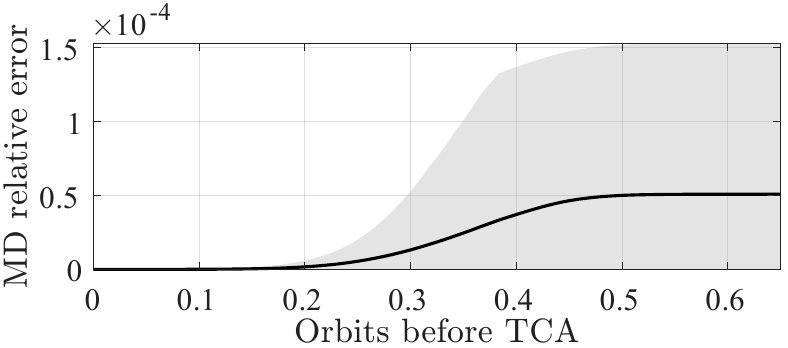}
    \caption{Statistics of the evolution of the relative \gls{md} error against a forward-propagation validation.}
    \label{fig:error_md}
\end{figure}

For each simulation, the control history computed with the \gls{gto} method is fed to a forward propagator, employing the same dynamics model (Keplerian) used to compute the polynomial maps, to validate the solution. The relative error of the metric computed with the Taylor polynomial with respect to the validated one is computed as
\begin{equation}
    e_i = \left|1 - \frac{\bar{d}^{(i)}_{(\cdot)}}{d_{(\cdot)}^{(val)}}\right| \ \ \ \forall \ t_i\in\{t_0, \ ..., \ t_N\}
\end{equation}
A Gaussian fit of this error is shown in \cref{fig:error_smd,fig:error_md} for the \gls{smd} and \gls{md} simulations, respectively. Across all simulations, $68\%$ of the solutions exhibit a relative error below $0.06\%$ in \gls{smd} and below $0.015\%$ ($30$~\si{cm}) in \gls{md}; with $99.7\%$ confidence, the errors are below $0.15\%$ in \gls{smd} and below $0.04\%$ in \gls{md}.

\begin{figure}
    \centering
    \includegraphics[width=0.9\linewidth]{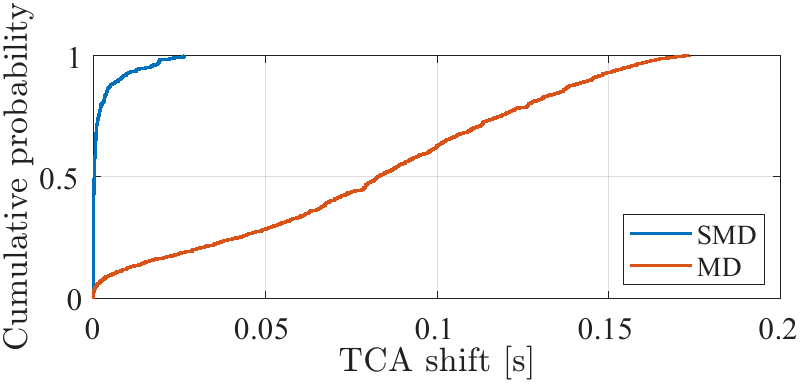}
    \vspace{-3mm}
    \caption{Cumulative probability of the \gls{tca} shift caused by the maneuver over the 2170 test cases.}
    \label{fig:tca_shift}
\end{figure}

In \cref{fig:tca_shift}, the cumulative distribution of the \gls{tca} offset caused by the maneuver is shown. The conjunction typically moves less for mostly radial maneuvers, while it changes faster for long tangential ones. Indeed, the \glspl{cam} computed with \gls{smd} show negligible change in \gls{tca} (always below $30$~\si{ms}); \gls{md} maneuvers, instead, which typically require longer thrusting windows, cause a maximum shift of around $170$~\si{ms}.

\begin{figure}
    \centering
    \includegraphics[width=0.9\linewidth]{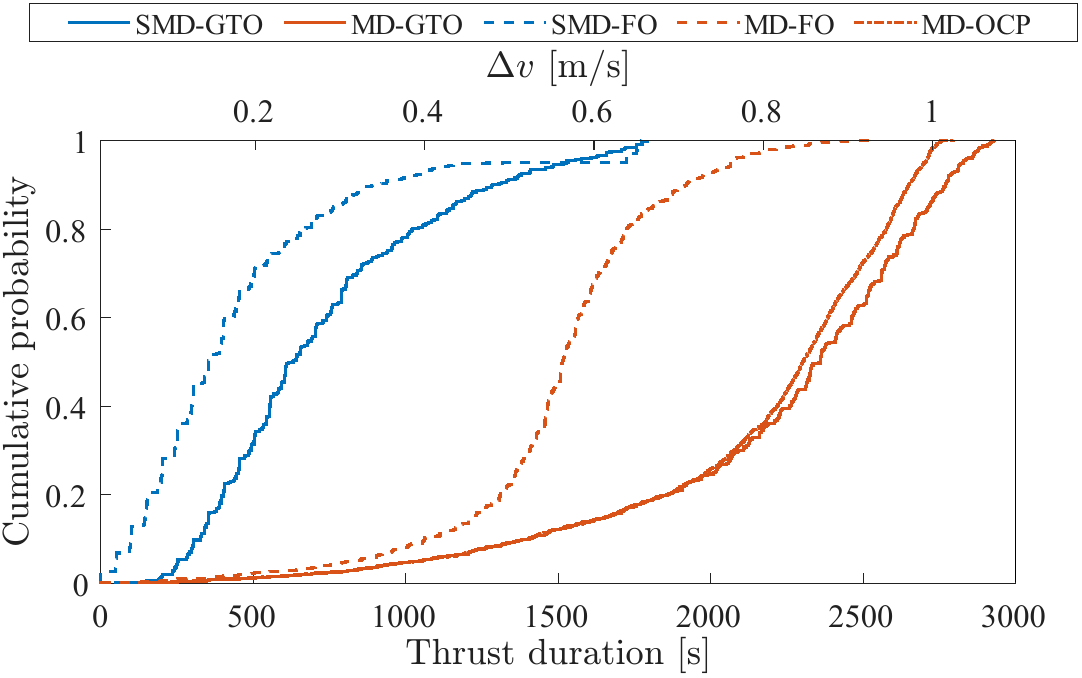}
    \caption{Cumulative probability of the maneuver duration over the 2170 test cases. The dashed lines indicate the \gls{fo} benchmark.}
    \label{fig:start_time}
\end{figure}

\begin{figure}
    \centering
    \includegraphics[width=0.9\linewidth]{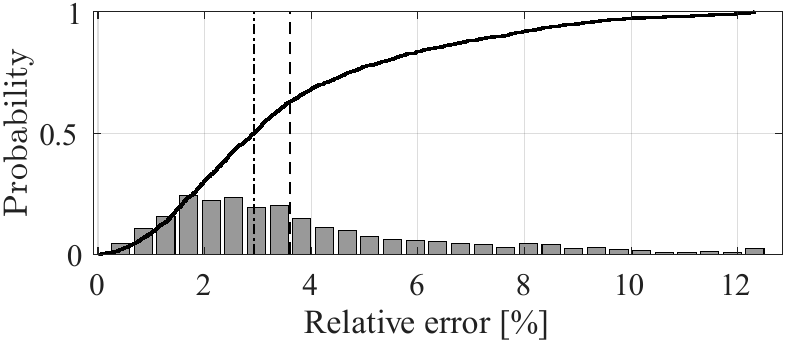}
    \caption{Probability distributions of the relative difference in $t_0$ results between the \gls{gto} and the \gls{ocp} methods. The median ($2.96\%$) and mean ($3.60\%$) are highlighted by vertical lines.}
    \label{fig:t0_OCP_comp_cdf}
\end{figure}

To assess the performance of the \gls{gto} method in terms of minimizing the thrust initiation time, the resulting $t_0$ values for the 2170 \gls{leo} conjunctions are compared to those using the \gls{ocp} method from Ref. \cite{dell2024characterization}. This method, based on optimal control theory, employs Pontryagin's Maximum Principle to find the true optimal solution. For the purpose of this work, a linearised version of this optimal method was used with near-perfect coincidence with the true optimal solution within the region of validity, which includes all 2170 test cases. The results of this comparison, which is only performed for the \gls{md} case, are shown in \cref{fig:start_time,fig:t0_OCP_comp_cdf}. In \cref{fig:start_time}, the \gls{md}-\gls{gto} curve and the \gls{md}-\gls{ocp} one are close to each other. The differences shown in \cref{fig:t0_OCP_comp_cdf} quantify the sub-optimality of the \gls{gto} method. The relative error is low for almost all conjunctions, with a median of $2.96\%$. In these figures, the $95^\text{th}$ percentiles are shwown, for which the error is less than $11\%$, though there are a number of outliers, with a maximum error at $30.5\%$. Nevertheless, this comparison shows that the \gls{gto} method forms a good approximation of the true optimal solution.

\begin{figure}
    \centering
    \includegraphics[width=\linewidth]{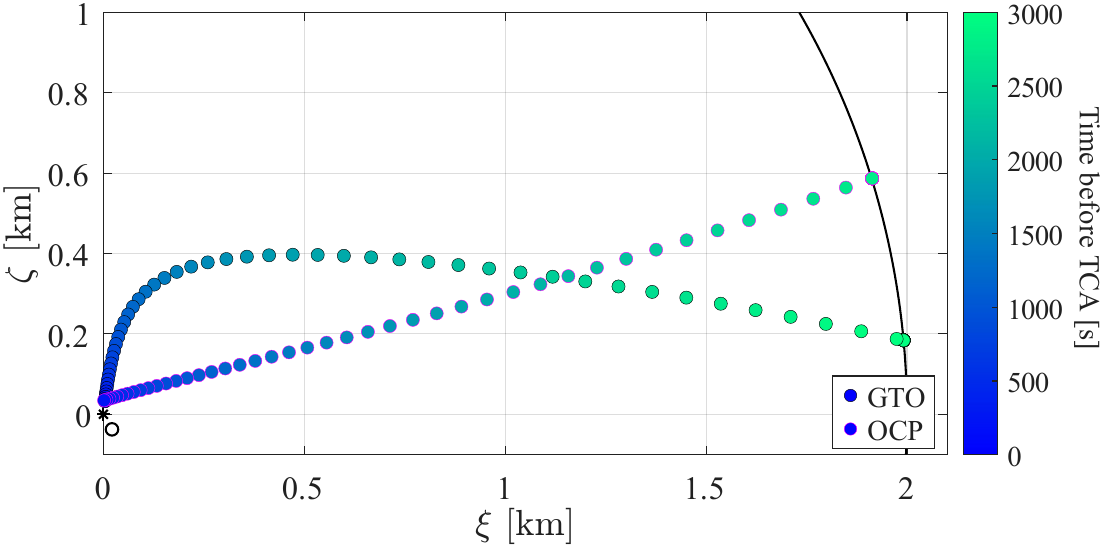}
    \vspace*{-5mm}
    \caption{Example (case 6/2170) of the projected evolution of the relative position in the B-plane during the maneuver.}
    \label{fig:bplane_oc}
\end{figure}

\cref{fig:bplane_oc} shows the different trajectories between the \gls{gto} and \gls{ocp} methods on the B-plane for a case where the difference in initiation time relatively largest within the $95^\text{th}$ percentile. In this case, the greedy trajectory direction differs significantly from the optimal solution, though both are valid \glspl{cam}. An analysis of the control profiles for this case indicates that both methods primarily rely on tangential thrust for the \gls{cam}, but the \gls{gto} method makes more use of normal thrust when close to \gls{tca}, while the \gls{ocp} method mostly uses radial thrust, which is more effective. 

\begin{figure}
    \centering
    \includegraphics[width=\linewidth]{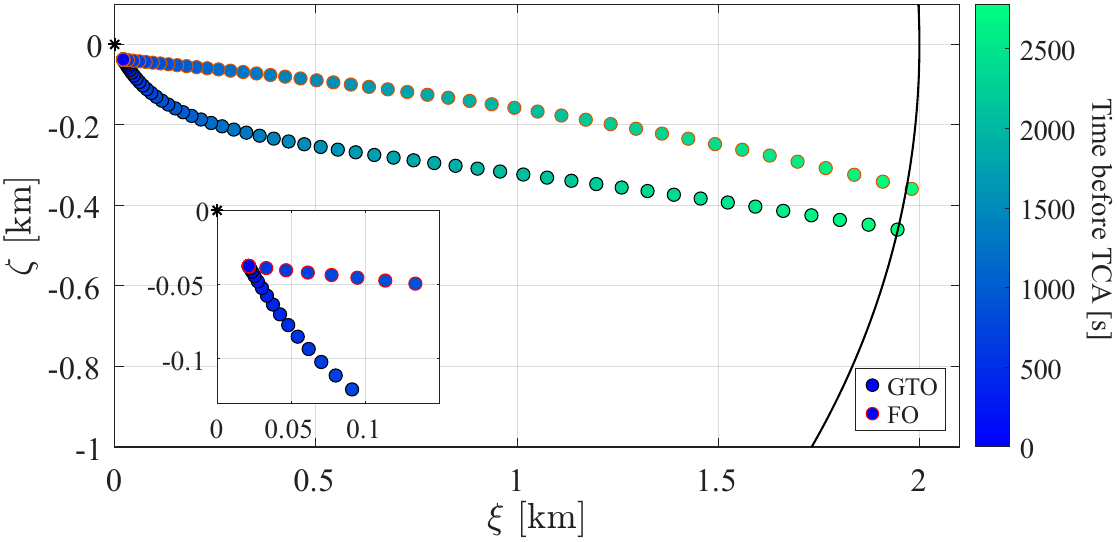}
    \vspace*{-5mm}
    \caption{Example (case 1/2170) of the projected evolution of the relative position in the B-plane during the maneuver.}
    \label{fig:bplane}
\end{figure}

\begin{figure}
    \centering
    \includegraphics[width=\linewidth]{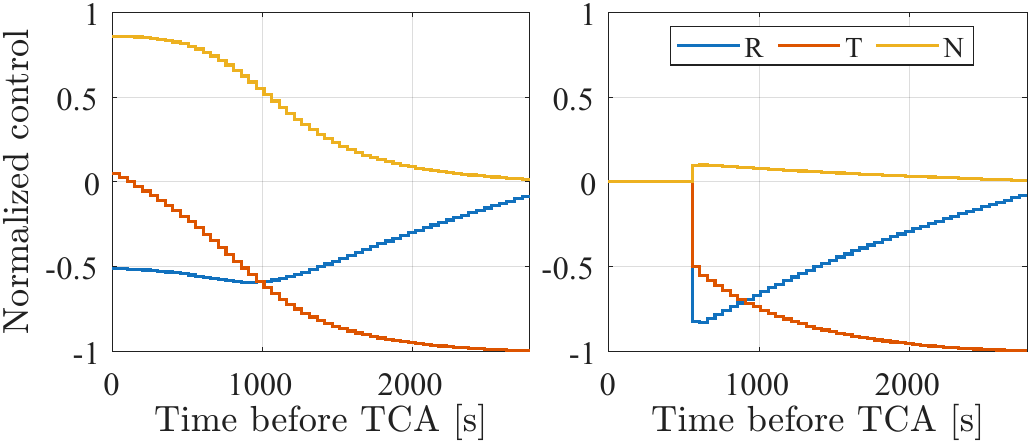}
    \vspace{-5mm}
    \caption{Control required to get the projections in \cref{fig:bplane}. Left: \gls{gto} method, right: \gls{fo}.}
    \label{fig:control}
\end{figure}

Even though the objective of the optimization is to minimize the thrust initiation time, the results are compared against a \gls{fo} solution computed over a time horizon equal to the minimum warning time using the \gls{scp} optimizer developed in \cite{pavanello2025thesis_chap3}. This comparison provides an indication of the fuel sub-optimality associated with the \gls{gto} approach. In particular, if the $\Delta v$ of the \gls{fo} solution is lower than that obtained with the minimum-time approach, this indicates that the thrust is not continuously active over the full warning interval. Consequently, within the minimum time given by the \gls{gto}, a better \gls{fo} solution exists. This appears evident in \cref{fig:start_time}, where, for both metrics, the total $\Delta v$ required by the maneuver and the thrust duration are compared (the two figures of merit are superimposed because in a low-thrust scenario without mass loss, the $\Delta v$ is proportional to the thrusting window). In both \gls{smd} and \gls{md} cases, the \gls{fo} solutions outperform the \gls{gto} method, requiring an average of $0.179$~\si{m/s} (\gls{smd}) and $0.570$~\si{m/s} (\gls{md}) against the $0.293$~\si{m/s} and $0.847$~\si{m/s} of the \gls{gto} method. This corresponds to an average reduction in thrust time by $41.3\%$ in the \gls{smd} case, and $33.6\%$ in the \gls{md} case. An example of the different trajectories achieved by the two methods is shown in \cref{fig:bplane}: using the \gls{md} as a collision metric, the \gls{fo} method follows a more straightforward path towards the edge. The \gls{gto} method wastes some thrust near \gls{tca}, when the effect is at its minimum. In fact, in \cref{fig:control}, one can see that the \gls{fo} solution takes advantage of a tangential thrust for slightly longer, and it almost completely avoids the out-of-plane thrust which is found by the \gls{gto} method. This is a paradigmatic example of why the \gls{gto} method is suboptimal from a \gls{fo} perspective, since most of the out-of-plane thrust is unnecessary. In fact, as can be observed in Fig. 7 of Ref. \cite{pavanello2025cammary}, the norm of the gradient of the metric with respect to the control action becomes negligible in the vicinity of \gls{tca} when compared with its magnitude earlier along the orbit. Consequently, enforcing a coast arc near \gls{tca} is a reasonable strategy for improving fuel efficiency, as control effort expended in this region has only a marginal impact on the optimization metric.
However, in many cases, which are not shown for the sake of conciseness, the \gls{gto} thrust aligns well with the \gls{fo} one, even though thrusting in proximity of the conjunction is almost always avoided by the \gls{fo} solver.

\begin{figure}
    \centering
    \includegraphics[width=0.9\linewidth]{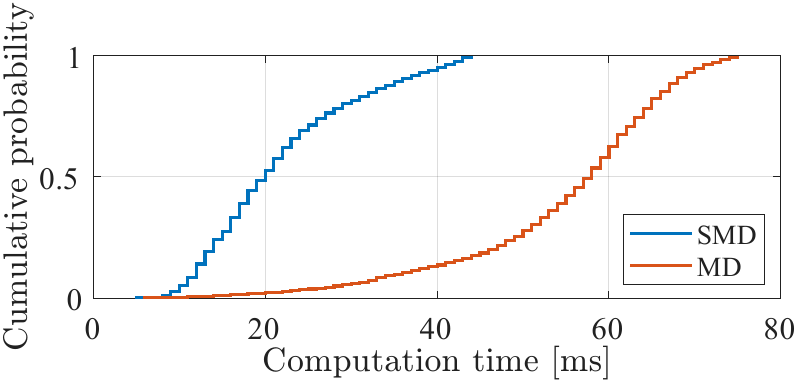}
    \caption{Cumulative probability of the computation time of the \gls{gto} method over the 2170 test cases.}
    \label{fig:comp_time}
\end{figure}

After this comparison, one might reconsider the need for the proposed \gls{gto} method. However, one should also be aware of the great computational efficiency of the \gls{gto} method: as shown in \cref{fig:comp_time}, runtimes are always kept below $80$~\si{ms}, with an average of $23$~\si{ms} and $56$~\si{ms} for the \gls{smd} and \gls{md} cases, respectively. This establishes the operational viability of the \gls{gto} method when warning time is a primary factor and fuel optimality can be sacrificed, considering the low absolute magnitude of the maneuvers. If time is not a factor, instead, operators could establish a pipeline in which a \gls{fo} method is executed after the \gls{gto}, using the minimum initiation time as input to the \gls{fo} problem. Moreover, the control profile solution from the \gls{gto} can serve as a first guess for \gls{fo} numerical solvers to speed up convergence, if needed.

\section{Conclusions}
\glsresetall
\label{sec:conclusions}

This work presented a methodology to determine the latest possible initiation time for a \gls{to} low-thrust \gls{cam} given a short warning time. The  \gls{gto} approach is based on a backward sweep in time, starting from the nominal \gls{tca} and iteratively stacking greedy control actions until the chosen safety threshold, expressed as either the \gls{smd} or the \gls{md} in the B-plane, is exceeded. Differential algebra is employed to construct polynomial maps of the system dynamics, enabling the efficient propagation of state sensitivities and the online tracking of the \gls{tca} shift induced by the maneuver, without re-propagating the dynamics at each sweep step.

The method was tested on ESA's Collision Avoidance Challenge dataset of 2170 \gls{leo} conjunctions, achieving accurate, computationally efficient results. Validation against a forward propagator showed that $99.7\%$ of solutions achieved a relative metric error below $0.15\%$. Computation times were kept below $80$~\si{ms} in all cases, with averages of $23$~\si{ms} and $56$~\si{ms} for the \gls{smd} and \gls{md} formulations, respectively, establishing the operational viability of the approach for on-board implementation. Comparison with an optimal control method showed a very limited median overestimation on the thrust time, while comparison with a fuel-optimal benchmark revealed that the \gls{gto} formulation requires, on average, between $33\%$ and $41\%$ more $\Delta v$. This sub-optimality is an acceptable trade-off when warning time is the primary constraint, and the \gls{fo} solution can further serve as an initial guess for \gls{fo} solvers.

\begin{footnotesize}
\bibliography{references}
\end{footnotesize}

\begin{IEEEbiography}[{\includegraphics[width=1in,height=1.25in,clip,keepaspectratio]{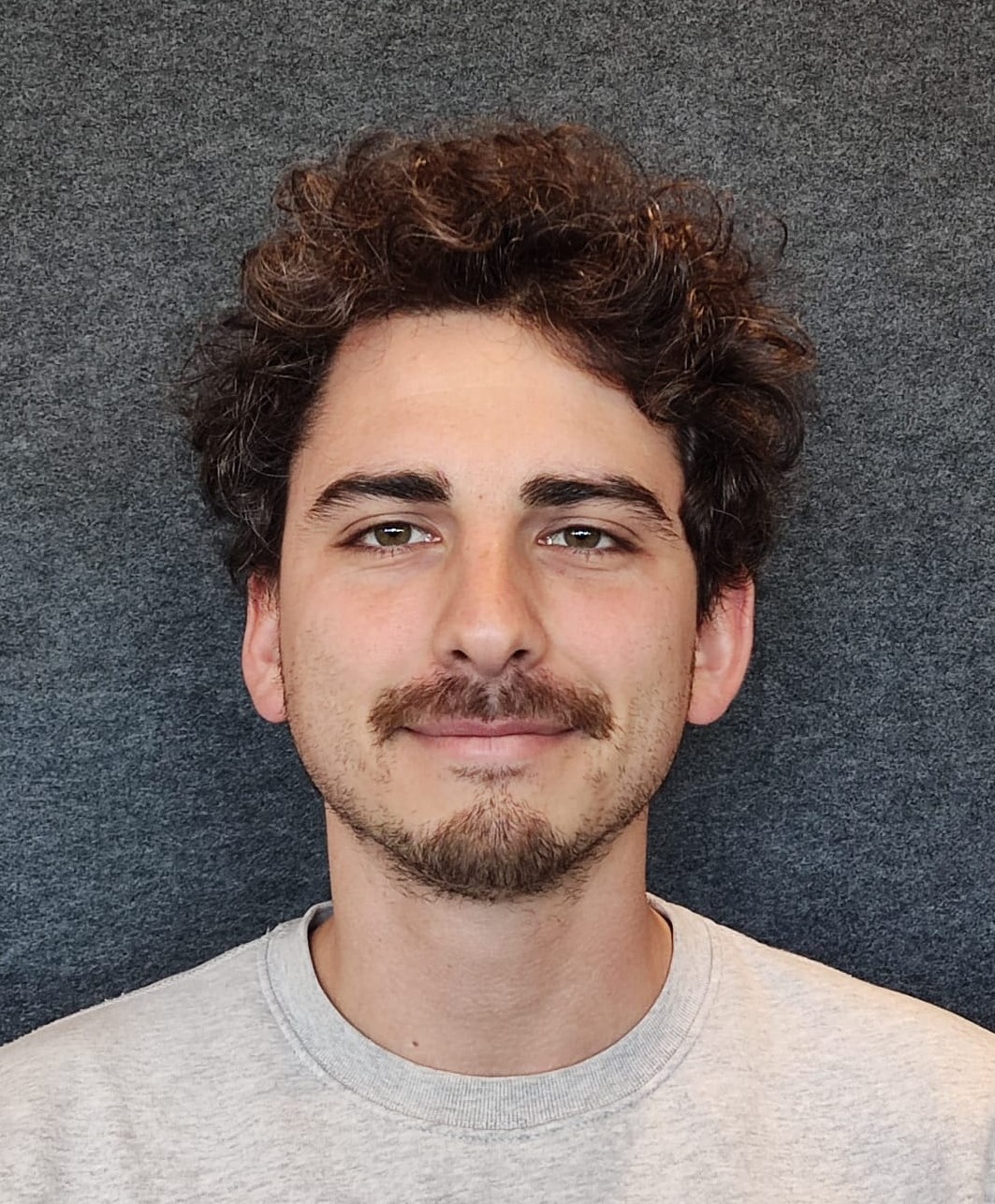}}]{Zeno Pavanello}{\space} obtained his MSc double degree in Aerospace Engineering in 2020 from the University of Padua and the Instituto Superior Técnico of Lisbon. He received his PhD degree from the University of Auckland in October 2025. He currently investigates space situational awareness and space traffic management-related topics as a PostDoc researcher at Politecnico di Milano, Italy.
\end{IEEEbiography}

\begin{IEEEbiography}
[{\includegraphics[width=1in,height=1.25in,clip,keepaspectratio]{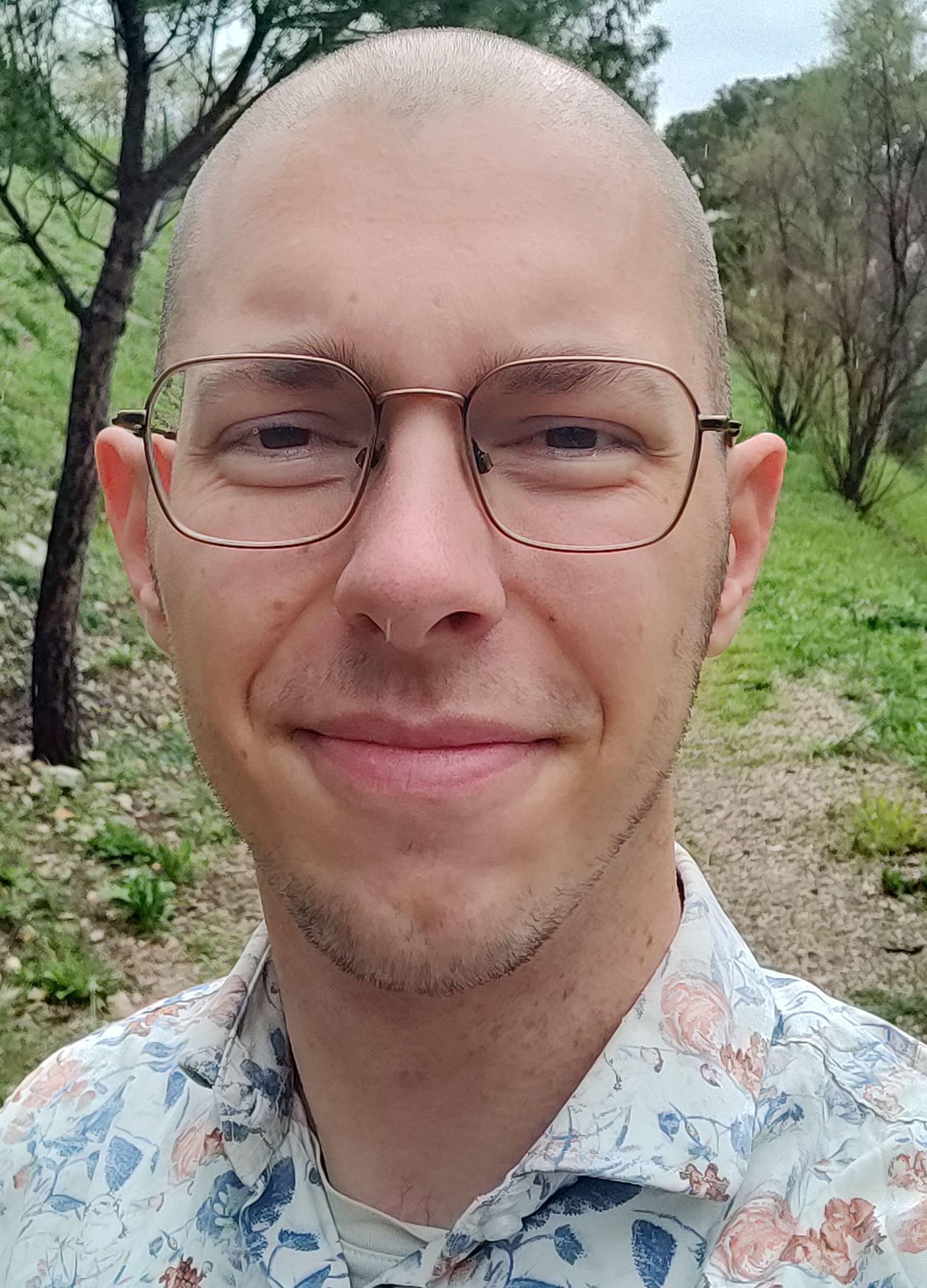}}]{Frank de Veld} received his MSc degree in 2022 at Delft University of Technology in aerospace engineering and his PhD degree on satellite control at Université Côte d'Azur in December 2025. Currently, he is active as a space situational awareness specialist for the Dutch Aerospace Research Center (NLR).
\end{IEEEbiography}

\begin{IEEEbiography}
[{\includegraphics[width=1in,height=1.25in,clip,keepaspectratio]{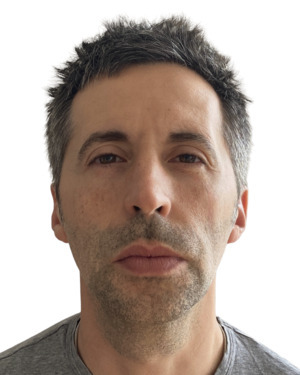}}]{Roberto Armellin} received his MSc and PhD degrees in aerospace engineering from Politecnico di Milano, Italy, in 2003 and 2007, respectively. Since November 2020, he has been a professor at Te P\=unaha \=Atea – Space Institute, the University of Auckland. His current research interests include space trajectory optimization, spacecraft navigation and guidance, and space situational awareness. 
\end{IEEEbiography}

\end{document}

%% file: Figures/back_sweep.tex
\tikzset{every picture/.style={line width=0.75pt}} 

\begin{tikzpicture}[x=0.75pt,y=0.75pt,yscale=-1,xscale=1]

\draw [color={rgb, 255:red, 181; green, 180; blue, 182 }  ,draw opacity=1 ] [dash pattern={on 0.84pt off 2.51pt}]  (200.2,96.12) .. controls (248.12,56.98) and (336.07,57.05) .. (362.07,64.33) ;
\draw [color={rgb, 255:red, 144; green, 19; blue, 254 }  ,draw opacity=1 ]   (83.2,260.6) .. controls (87.2,222.57) and (116.2,177.77) .. (123.2,171.77) ;
\draw [color={rgb, 255:red, 181; green, 180; blue, 182 }  ,draw opacity=1 ] [dash pattern={on 0.84pt off 2.51pt}]  (123.2,171.77) .. controls (151.2,125.48) and (185.2,100.48) .. (201.2,94.12) ;
\draw [color={rgb, 255:red, 144; green, 19; blue, 254 }  ,draw opacity=1 ]   (133.2,179.6) .. controls (155.2,143.48) and (204.01,105.48) .. (209.51,104.37) ;
\draw [color={rgb, 255:red, 144; green, 19; blue, 254 }  ,draw opacity=1 ]   (218.2,121.6) .. controls (243.2,97.48) and (319.12,85.98) .. (345.2,90.87) ;
\draw [color={rgb, 255:red, 181; green, 180; blue, 182 }  ,draw opacity=1 ] [dash pattern={on 0.84pt off 2.51pt}]  (207.7,104.37) .. controls (249.12,72.98) and (325.07,71.33) .. (354.07,75.33) ;
\draw    (83.2,260.6) .. controls (108.2,185.6) and (206.2,81.6) .. (355.2,106.6) ;
\draw  [fill={rgb, 255:red, 0; green, 0; blue, 0 }  ,fill opacity=1 ] (121.39,171.77) .. controls (121.39,172.77) and (122.2,173.58) .. (123.2,173.58) .. controls (124.2,173.58) and (125.01,172.77) .. (125.01,171.77) .. controls (125.01,170.77) and (124.2,169.96) .. (123.2,169.96) .. controls (122.2,169.96) and (121.39,170.77) .. (121.39,171.77) -- cycle ;
\draw [color={rgb, 255:red, 15; green, 95; blue, 190 }  ,draw opacity=1 ][line width=1.5]    (77.2,258.6) -- (88.2,262.6) ;
\draw [color={rgb, 255:red, 15; green, 95; blue, 190 }  ,draw opacity=1 ][line width=1.5]    (129.2,176.6) -- (138.2,183.87) ;
\draw [color={rgb, 255:red, 15; green, 95; blue, 190 }  ,draw opacity=1 ][line width=1.5]    (215.2,116.4) -- (221.2,125.87) ;
\draw [color={rgb, 255:red, 15; green, 95; blue, 190 }  ,draw opacity=1 ][line width=1.5]    (338.2,97.38) -- (336.2,109.87) ;
\draw [color={rgb, 255:red, 0; green, 0; blue, 0 }  ,draw opacity=1 ]   (138.2,183.87) -- (123.2,171.77) ;
\draw [color={rgb, 255:red, 0; green, 0; blue, 0 }  ,draw opacity=1 ]   (201.2,94.12) -- (221.2,125.87) ;
\draw [color={rgb, 255:red, 0; green, 0; blue, 0 }  ,draw opacity=1 ]   (345.2,90.87) -- (343.2,104.87) ;
\draw  [fill={rgb, 255:red, 0; green, 0; blue, 0 }  ,fill opacity=1 ] (205.89,104.37) .. controls (205.89,105.37) and (206.7,106.18) .. (207.7,106.18) .. controls (208.7,106.18) and (209.51,105.37) .. (209.51,104.37) .. controls (209.51,103.37) and (208.7,102.56) .. (207.7,102.56) .. controls (206.7,102.56) and (205.89,103.37) .. (205.89,104.37) -- cycle ;
\draw  [fill={rgb, 255:red, 0; green, 0; blue, 0 }  ,fill opacity=1 ] (343.39,90.87) .. controls (343.39,91.87) and (344.2,92.68) .. (345.2,92.68) .. controls (346.2,92.68) and (347.01,91.87) .. (347.01,90.87) .. controls (347.01,89.87) and (346.2,89.06) .. (345.2,89.06) .. controls (344.2,89.06) and (343.39,89.87) .. (343.39,90.87) -- cycle ;
\draw  [fill={rgb, 255:red, 0; green, 0; blue, 0 }  ,fill opacity=1 ] (199.39,94.31) .. controls (199.39,95.31) and (200.2,96.12) .. (201.2,96.12) .. controls (202.2,96.12) and (203.01,95.31) .. (203.01,94.31) .. controls (203.01,93.31) and (202.2,92.5) .. (201.2,92.5) .. controls (200.2,92.5) and (199.39,93.31) .. (199.39,94.31) -- cycle ;
\draw  [fill={rgb, 255:red, 0; green, 0; blue, 0 }  ,fill opacity=1 ] (352.26,74.52) .. controls (352.26,75.52) and (353.07,76.33) .. (354.07,76.33) .. controls (355.07,76.33) and (355.88,75.52) .. (355.88,74.52) .. controls (355.88,73.53) and (355.07,72.72) .. (354.07,72.72) .. controls (353.07,72.72) and (352.26,73.53) .. (352.26,74.52) -- cycle ;
\draw  [fill={rgb, 255:red, 0; green, 0; blue, 0 }  ,fill opacity=1 ] (360.26,64.33) .. controls (360.26,65.33) and (361.07,66.14) .. (362.07,66.14) .. controls (363.07,66.14) and (363.88,65.33) .. (363.88,64.33) .. controls (363.88,63.33) and (363.07,62.52) .. (362.07,62.52) .. controls (361.07,62.52) and (360.26,63.33) .. (360.26,64.33) -- cycle ;
\draw [color={rgb, 255:red, 0; green, 0; blue, 0 }  ,draw opacity=1 ]   (354.07,75.33) -- (349.07,105.33) ;
\draw [color={rgb, 255:red, 0; green, 0; blue, 0 }  ,draw opacity=1 ]   (362.07,64.33) -- (355.2,106.6) ;

\draw (71.33,264.67) node [anchor=north west][inner sep=0.75pt]    {$\vec{x}_{N-3}$};
\draw (89.33,251.67) node [anchor=north west][inner sep=0.75pt]  [color={rgb, 255:red, 15; green, 95; blue, 190 }  ,opacity=1 ]  {$t_{N-3}$};
\draw (142.33,174.67) node [anchor=north west][inner sep=0.75pt]  [color={rgb, 255:red, 15; green, 95; blue, 190 }  ,opacity=1 ]  {$t_{N-2}$};
\draw (220.2,121.6) node [anchor=north west][inner sep=0.75pt]  [color={rgb, 255:red, 15; green, 95; blue, 190 }  ,opacity=1 ]  {$t_{N-1}$};
\draw (325.33,109.67) node [anchor=north west][inner sep=0.75pt]  [color={rgb, 255:red, 15; green, 95; blue, 190 }  ,opacity=1 ]  {$t_{N}$};
\draw (111.33,212.67) node [anchor=north west][inner sep=0.75pt]  [color={rgb, 255:red, 144; green, 19; blue, 254 }  ,opacity=1 ]  {$\vec{u}_{N-3}$};
\draw (183.33,144.67) node [anchor=north west][inner sep=0.75pt]  [color={rgb, 255:red, 144; green, 19; blue, 254 }  ,opacity=1 ]  {$\vec{u}_{N-2}$};
\draw (273.33,111.67) node [anchor=north west][inner sep=0.75pt]  [color={rgb, 255:red, 144; green, 19; blue, 254 }  ,opacity=1 ]  {$\vec{u}_{N-1}$};
\draw (88.33,153.67) node [anchor=north west][inner sep=0.75pt]    {$\delta \vec{x}_{N-2}$};
\draw (166.33,75.67) node [anchor=north west][inner sep=0.75pt]    {$\delta \vec{x}_{N-1}$};
\draw (347.33,46.67) node [anchor=north west][inner sep=0.75pt]    {$\delta \vec{x}_{N}$};

\end{tikzpicture}